\documentclass{amsart}

\usepackage[frame,cmtip,arrow,matrix,line,graph,curve]{xy}
\usepackage{amssymb}
\usepackage{latexsym}
\usepackage{euscript}

\newtheorem{theorem}{Theorem}[section]
\newtheorem{proposition}[theorem]{Proposition}
\newtheorem{corollary}[theorem]{Corollary}
\newtheorem{lemma}[theorem]{Lemma}
\newtheorem{conjecture}[theorem]{Conjecture}

\theoremstyle{definition}
\newtheorem{definition}[theorem]{Definition}
\newtheorem{example}[theorem]{Example}
\newtheorem{assumption}[theorem]{Assumption}

\theoremstyle{remark}
\newtheorem{remark}[theorem]{Remark}

\newcommand{\pp}{\mathbb{P}}
\newcommand{\qq}{\mathbb{Q}}
\newcommand{\cc}{\mathbb{C}}
\newcommand{\rr}{\mathbb{R}}
\newcommand{\zz}{\mathbb{Z}}

\newcommand{\lik}{\mathfrak{k}}

\newcommand{\lih}{\mathfrak{h}}
\newcommand{\liu}{\mathfrak{u}}
\newcommand{\lit}{\mathfrak{t}}

\newcommand{\id}{\textrm{id}}
\newcommand{\Hom}{\mathrm{Hom}}
\newcommand{\Lie}{\mathrm{Lie \,}}

\newcommand{\codim}{\mathrm{codim \,}}
\newcommand{\Stab}{\mathrm{Stab} \,}

\newcommand{\git}{/\!\!/}
\newcommand{\udot}{\,{\bf \dot{}}\,}

\newcommand{\ck}[1]{\mathcal{C}\udot_{\!\! K}(#1)}
\newcommand{\eqc}[2]{\mathcal{C}\udot_{\!\! #1}(#2)}

\newcommand{\ic}[1]{\mathcal{IC}\udot (#1)}
\newcommand{\ict}[1]{\mathcal{IC}\udot_{\!\! t} (#1)}
\newcommand{\icp}[2]{\mathcal{IC}\udot_{\!\! #1} (#2)}
\newcommand{\con}[1]{\mathrm{Con}\udot (#1)}

\newcommand{\dgs}[1]{\mathcal{#1}\udot}
\newcommand{\dgsx}[2]{\mathcal{#1}\udot(#2)}
\newcommand{\dgss}[2]{\mathcal{#1}\udot_{\!\! #2}}
\newcommand{\ve}[1]{\dgs{V}(#1)}

\newcommand{\supp}{\mathrm{supp} \,}

\newcommand{\rdf}[1]{R#1}

\newcommand{\der}[1]{{\bf D}(#1)}
\newcommand{\pdera}[3]{{\bf D}^{#1}_{\geq#2}(#3)}
\newcommand{\pderb}[3]{{\bf D}^{#1}_{\leq#2}(#3)}

\newcommand{\ta}[1]{\tau_{>#1}}
\newcommand{\tgeq}[1]{\tau_{\geq #1}}
\newcommand{\tbon}[2]{\tau_{<#1}^{#2}}
\newcommand{\tb}[1]{\tau_{<#1}}
\newcommand{\ptb}[2]{\tau^{#1}_{\leq#2}}
\newcommand{\pta}[2]{\tau^{#1}_{\geq#2}}

\newcommand{\hyp}[3]{\mathbb{H}^{#1}(#2 ; #3)}

\newcommand{\hk}[2]{H^{#1}_K(#2)}
\newcommand{\ih}[2]{I\!H^{#1}(#2)}
\newcommand{\iht}[2]{I\!H_t^{#1}(#2)}
\newcommand{\ihp}[3]{I\!H_{#1}^{#2}(#3)}
\newcommand{\h}[2]{H^{#1}(#2)}
\newcommand{\sh}[2]{H^{#1}(#2)}

\begin{document}
\pagestyle{plain}

\title{The cosupport axiom, equivariant cohomology and the intersection cohomology of certain symplectic quotients}
\date{September 2000}
\author{Young-Hoon Kiem and Jonathan
Woolf}
\address{Young-Hoon Kiem, Stanford
University, CA 94305, USA; kiem@math.stanford.edu} 
\address{Jonathan Woolf, Christ's College, Cambridge, CB2 3BU,
UK; jw301@cus.cam.ac.uk}
\thanks{Y.H. Kiem was partially supported by
an Alfred P. Sloan doctoral dissertation fellowship.}
\thanks{J. Woolf was supported by the William G. Seggie-Brown research fellowship.}
\maketitle
\section*{Introduction}
Suppose $K\times M \to M$ is a Hamiltonian action of a compact Lie group $K$ on a symplectic manifold $M$, and that $\mu:M \to \lik^*$ is a proper equivariant moment map for this action. To avoid unnecessary symmetry we assume throughout this paper that there is at least one point in $\mu^{-1}(0)$ whose stabilizer is zero dimensional. The symplectic quotient, or reduction at zero, is given by the construction $M_0 := Z/K$ where $Z=\mu^{-1}(0)$. The infinitesimal orbit type stratification of $Z$ descends to a stratification of $M_0$ by symplectic orbifolds. When $0$ is a regular value of the moment map there is a single stratum. In this case the cohomology of $M_0$ can be studied via the surjective restriction:
$$
\hk{*}{M} \longrightarrow \hk{*}{Z} \cong  \h{*}{M_0}.
$$
When $0$ is no longer a regular value then, in general, $\hk{*}{Z}$ will have infinite dimension, and the pull-back $\h{*}{M_0} \to \hk{*}{Z}$ will be neither surjective nor injective. 

In this paper we show how, under certain conditions on the action, we can still extract topological information about $M_0$ from $\hk{*}{Z}$ when the reduction is no longer regular. The information we obtain is a complete description of the middle perversity intersection cohomology $\ih{*}{M_0}$ as a graded vector space with non-degenerate pairing between classes of complementary degree. In many cases when dealing with singular spaces this appears to be more natural data  than the ordinary cohomology. 

Intersection cohomology has its most elegant description in terms of the derived category of constructible sheaves. We are interested in the relationship between two objects, the derived push-down $\ck{Z}$ of the constant sheaf on $EK\times_KZ$, and the middle perversity intersection cohomology complex $\ic{M_0}$. The hypercohomology functor takes these respectively to $\hk{*}{Z}$ and $\ih{*}{M_0}$. In \S\ref{kirwan map} we show how the partial resolution constructed in \cite{mes} allows us to define (non-canonically) a morphism 
$$
\kappa : \ck{Z} \longrightarrow \ic{M_0}.
$$
This generalises a construction first used by Kirwan in \cite{k3}. 
 
The intersection cohomology complex $\ic{M_0}$ is characterised as an object of the derived category by three axioms which we call normalisation, support and cosupport. The normalisation axiom, that a complex is quasi-isomorphic to the constant sheaf on the unique open stratum, is always satisfied by $\ck{Z}$. A simple check shows that the support axiom holds for $\ck{Z}$ precisely when $K$ has only finite stabilizers on $Z$, or, equivalently, the reduction is regular. (Of course cosupport also holds for $\ck{Z}$ in this situation because it is quasi-isomorphic to $\ic{M_0}$ under the regularity assumption.) In this paper we demand rather less, we allow support to fail, so the reduction need not be regular, but demand that cosupport holds for $\ck{Z}$. We say an action for which this occurs is \emph{almost-balanced}. The reason for the terminology is that the almost-balanced condition is equivalent to asking that the weights of the action on a normal fibre to an infinitesimal orbit type stratum be evenly distributed about the origin. This is explained more precisely in \cite{kiem}. As well as their theoretical interest there are interesting examples of almost-balanced actions. For instance the moduli space of bundles over a Riemann surface arises as a quotient by an almost-balanced action. The rank $2$ case is discussed in some detail in \cite{kiem2}. The `standard' example of $SU_2$ acting on $(\pp^1)^n$ is also almost-balanced and we discuss this in \S \ref{simple example}. Very loosely, one expects actions arising in highly symmetric contexts to be almost-balanced. 

We show in \S\ref{cosupport condition} that when the action is almost-balanced there is a canonical morphism 
$$
\lambda : \ic{M_0} \longrightarrow \ck{Z}.
$$
which splits $\kappa$ i.e. $\kappa\lambda =\id$. If we make the slightly stronger assumption that the action is weakly balanced we can identify the subspace of $\hk{*}{Z}$ which corresponds to $\ih{*}{M_0}$. The naturality of $\lambda$ ensures that the intersection pairing corresponds to the normal product on equivariant cohomology. 

Finally, in \S\ref{circle actions}, we show that the same set of ideas can be applied to \emph{any} circle action. The key point is that we can define a non-middle perversity $n$ for which $\ihp{n}{*}{M_0} \cong \ih{*}{M_0}$. The statements of the previous section then hold for arbitrary circle actions but with $n$ replacing the middle perversity.     
\tableofcontents
\section{Intersection cohomology}
\label{intersection}
Suppose $X$ is a connected pseudomanifold of dimension $2n$ with a fixed even dimensional topological stratification $\{S_\alpha\}$. Let $\der{X}$ be the bounded below derived category of constructible sheaves on $X$. This is a triangulated category with shift functor $[n]$ and truncation functors $\tb{n}$ and $\ta{n}$. Slightly more general than $\tb{n}$ is the functor given by truncation on a closed subset, which was introduced by Deligne and defined in \cite[\S1.14]{gm2}. Let $Y$ be a closed subset of $X$. Then there is a functor $\tbon{n}{Y}$ and a morphism  $\tbon{n}{Y}\dgs{A} \to \dgs{A}$ which induces
\begin{equation}
\label{truncation}
H^i_x(\tbon{n}{Y}\dgs{A}) \cong \left\{
\begin{array}{ll}
0 & x\in Y \textrm{ and } i \geq n;\\
H^i_x(\dgs{A}) & \textrm{otherwise}
\end{array}
\right.
\end{equation}
for any sheaf complex $\dgs{A}$. Of course $\tbon{n}{X} = \tb{n}$. 

We denote the cone on a morphism $\alpha$ in $\der{X}$ by $\con{\alpha}$. The category $\der{X}$ also has a tensor product and Verdier duality functor. Stratified maps between spaces induce the standard derived functors $f^*,\rdf{f_*},f^!$ and $\rdf{f_!}$ of sheaf theory.  A reasonably complete description of this category is given in \cite[\S1]{gm2}.
\begin{definition}
A \emph{perversity} $p$ is a function $\{S_\alpha\} \to \zz$. We say $p$ is \emph{monotone} if $p(S_\alpha) \leq p(S_\beta)$ whenever $\overline{S_\alpha} \supset S_\beta$. 
\end{definition}
Given a perversity $p$ we define full subcategories of $\der{X}$ by
\begin{eqnarray*}
\pderb{p}{0}{X} & = & \{ \dgs{A} \ |\ \sh{i}{\jmath_\alpha^*\dgs{A}} =0 \textrm{ for } i > p(S_\alpha)\} \\
\textrm{and }\pdera{p}{0}{X} & = & \{ \dgs{A} \ |\ \sh{i}{\jmath_\alpha^!\dgs{A}} =0 \textrm{ for } i < p(S_\alpha)\}.
\end{eqnarray*}
\begin{theorem}[{\cite[Ch. 7, 1.2.1]{kos}}]
The pair $(\pderb{p}{0}{X},\pdera{p}{0}{X})$ of subcategories defines a $\bold{t}$-structure on $\der{X}$. There are left and right adjoints $\ptb{p}{0}$ and $\pta{p}{0}$ respectively to the inclusions
$$\pderb{p}{0}{X} \hookrightarrow \der{X} \quad \textrm{ and } \quad \pdera{p}{0}{X} \hookrightarrow \der{X}.
$$
\end{theorem}
Now suppose $p$ is a monotone perversity with $p(U)=0$ where $U$ is the unique open stratum. We define two further perversities $p^+$ and $p^-$ by 
$$
p^\pm(S_\alpha) = \left\{ 
\begin{array}{ll}
0 & S_\alpha = U \\
p(S_\alpha) \pm 1 & \textrm{otherwise.}
\end{array}
\right.
$$
\begin{theorem}[{\cite[\S 3]{gm2}}]
\label{ic exists}
There is a unique object 
$$
\icp{p}{X} \in \pderb{p^-}{0}{X} \cap \pdera{p^+}{0}{X}
$$ 
with the property that $\icp{p}{X}|_U$ is quasi-isomorphic to the constant sheaf $\qq$ on $U$. 
\end{theorem}
We call $\icp{p}{X}$ the perversity $p$ intersection cohomology complex. Its hypercohomology $\ihp{p}{*}{X}$ is the perversity $p$ intersection cohomology.
The uniqueness part of the proof rests on the following lemma which is proved by induction over the stratification:
\begin{lemma}
\label{unique morphisms} Suppose $\dgs{A}|_U$ is quasi-isomorphic to the constant sheaf $\qq$.
Then we have the following implications:
\begin{eqnarray*}
\dgs{A} \in \pderb{p^-}{0}{X} & \implies & \Hom(\dgs{A},\icp{p}{X}) \cong \qq \\\dgs{A} \in \pderb{p^+}{0}{X} & \implies & \Hom(\icp{p}{X},\dgs{A}) \cong \qq. 
\end{eqnarray*}
\end{lemma}
We define the \emph{middle perversity} to be the function $m(S_\alpha) = \frac{1}{2}\codim S_\alpha$. We will denote the  middle perversity intersection cohomology complex simply by $\ic{X}$. The \emph{top perversity} is given by $t(S_\alpha)=2m(S_\alpha)$. A simple check shows that 
$\ic{X}^{\otimes 2} \in \pderb{t^-}{0}{X}$
and so lemma \ref{unique morphisms} guarantees there is a natural morphism to $\ict{X}$. 
\begin{theorem}[{\cite[\S 5.2]{gm2}}]
This morphism gives rise to pairings
\[
\xymatrix{
\ih{n-i}{X} \otimes \ih{n+i}{X} \ar[r]& \iht{2n}{X} \ar[r]^{\quad \int_X} & \qq
}
\]
which are non-degenerate when $X$ is compact.
\end{theorem}
To simplify the notation we say that 
\begin{eqnarray*}
\dgs{A} \textrm{ satisfies \emph{support} } & \iff & \dgs{A} \in \pderb{m^-}{0}{X} \\
\dgs{A} \textrm{ satisfies \emph{cosupport} } & \iff & \dgs{A} \in \pdera{m^+}{0}{X}.
\end{eqnarray*}
Thus theorem \ref{ic exists} says that $\ic{X}$ is characterised by the conditions that it satisfies support, cosupport and is the constant sheaf on $U$.
\section{Singular symplectic quotients}
\label{singular}
Suppose $(M,\omega)$ is a Hamiltonian $K$-space for some compact connected Lie
group $K$. By this we mean that $K$ acts on $M$ preserving the
symplectic form $\omega$ in such a way that there is an equivariant moment map 
$\mu : M \to \lik^*$ satisfying
$$
\langle d\mu , a \rangle = \imath_{a_M}\omega
$$
where $a \in \lik$ and $a_M$ is the vector field on $M$ arising from
the infinitesimal action of $a$. Throughout this paper we will assume that $\mu$ is proper, and in this
case call $M$ a \emph{proper Hamiltonian $K$-space}.

It is well known that if $0$ is a regular value of $\mu$ then the
topological space $M_0 = \mu^{-1}(0)/K$ can naturally be given the structure
of a compact symplectic orbifold. When $0$ is not a regular value Lerman and Sjamaar show in \cite{sl} that $M_0$ is a stratified
symplectic space. We briefly describe the
stratification but ignore the Poisson structure on the functions since
this is unnecessary for our topological applications.

Let us put $Z= \mu^{-1}(0)$. Suppose $H$ is a compact subgroup of $K$ with Lie algebra $\lih$. Let $Z_H$ be the subset of points of $Z$ whose stabilizer is precisely $H$. We also define subsets of $Z$
$$
Z_{\lih} = \{z \in Z : \Lie \Stab z = \lih\} \quad\textrm{and}\quad Z_{(\lih)} = \{z \in Z : \Lie \Stab z \in (\lih)\}
$$
where $(\lih)$ is the set of 
subalgebras of $\lik$ conjugate to $\lih$. 

For each nonempty subset $Z_{(\lih)}$ of $Z$ fix once and for all a representative $\lih$ of the conjugacy class $(\lih)$ and let $H$ be the generic stabiliser of a point in $Z_{\lih}$ (so that $\Lie(H) = \lih$). The set $I_Z$ of these Lie subalgebras forms an indexing set for a topological stratification $\mathcal{S}_\lik$ of $M_0$ by infinitesimal orbit types. The  strata are the sets consisting of those points representing orbits in the subsets
$Z_{(\lih)}$ for $\lih \in I_Z$. These strata are symplectic orbifolds, with a
symplectic structure induced from $\omega$. This stratification
makes $M_0$ into a compact topological pseudomanifold --- see \cite[\S3.1.1]{mes} for details. 
\subsection{Relation with algebraic quotients}
\label{git}
There is a well known principle that symplectic and algebraic
quotients are `the same'. Suppose $V$ is a smooth complex projective variety with a given
embedding $V \hookrightarrow \pp^n$.
Let $G$ be a reductive algebraic group which acts on $V$ via  a
homomorphism $\rho : G \rightarrow GL_{n+1}$. Geometric invariant
theory constructs a categorical quotient $V^{ss} \git G$ of the Zariski open subset $V^{ss}$ of semistable points. 

Let $K$ be a maximal compact subgroup of $G$. We may assume, by conjugating if necessary, that $K$
maps into $U_{n+1}$ under the homomorphism $\rho$. The Fubini-Study
form on $\pp^n$ restricts to a K\"ahler form on $V$ which is preserved
by the action of $K$. Further the action of $K$ is Hamiltonian and
there is a moment map
$$
\mu : V \hookrightarrow \pp^n \stackrel{\phi}{\longrightarrow} \liu_{n+1}^*
\stackrel{(d\rho)^*}{\longrightarrow} \lik^*
$$
where $\phi$ is the moment map for the standard action of $U_{n+1}$ on
$\pp^n$.
\begin{theorem}(\cite{k1,ne})
\label{global algebraicity}
The zero set $\mu^{-1}(0)$ of the moment map is contained within the
semistable points $V^{ss}$. This inclusion induces a homeomorphism
$$
\mu^{-1}(0) / K \longrightarrow V^{ss} \git G.
$$
\end{theorem}
More generally let $M$ be a proper Hamiltonian $K$-space. We fix a
choice of compatible $K$-invariant almost complex structure $J$ and
metric $g$ on $M$. Let $p \in Z_H \subset Z_{(\lih)}$. Define
$$
V_p = \big(T_p(Kp) \oplus JT_p(Kp)\big)^\perp \leq T_pM
$$
where $\perp$ denotes the orthogonal complement with respect to $g$,
and set
\begin{equation}
\label{w space}
W_p = (V_p)_\lih^\perp \leq V_p
\end{equation}
where $(V_p)_\lih$ is the subspace invariant under the infinitesimal
$\lih$ action. $W_p$ has a Hermitian
structure, induced from the triple $(J,g,\omega)$, with respect to which
it becomes a Hamiltonian $H$-space. It follows from the local normal
form for the moment map (see \cite{ma,gs1}) that a
neighbourhood of the point $q$ in $M_0$ representing the orbit
$Kp$ is homeomorphic to the product of the reduction at zero of $W_p$ with $\cc^r$ for some $r$. By a
small extension of theorem \ref{global algebraicity} to the quasi-projective
variety $W_p$ we deduce that a neighbourhood of $q$ in $M_0$ is
homeomorphic to the product of the geometric invariant theory quotient $W_p \git H^\cc$ with $\cc^r$. Thus all symplectic reductions are locally homeomorphic to algebraic varieties. 
\subsection{Partial desingularisation}
\label{partial}
In \cite{k2} Kirwan showed how a singular geometric invariant theory
quotient may be canonically desingularised by a particular sequence of
blowups. This approach was extended to singular reductions of symplectic
manifolds with proper moment map in \cite{mes}. Lack of
canonicity in the definition of \emph{symplectic} blowup means the
resulting desingularisation is certainly not unique up to
symplectomorphism. However it is unique up to homeomorphism and as we
are interested in purely topological questions in this paper this is
quite sufficient for our purposes.

Suppose that $M$ is a symplectic manifold with a Hamiltonian action of
a compact connected Lie group $K$ for which the moment map
$\mu:M \to \lik^*$ is proper. Let $\lih \in I_Z$ be a Lie subalgebra of $\lik$
indexing a stratum $S \in \mathcal{S}_\lik$ of maximal depth.  
\begin{proposition}[partial desingularisation]
\label{pdprop}
We can find a continuous surjection $\widetilde M_0 \stackrel{\pi}{\longrightarrow} M_0$ where $\widetilde M_0$ is the reduction of a proper Hamiltonian
$K$-space $\widetilde M$ such that
\begin{enumerate}
\item the restriction $\widetilde M_0 \setminus \pi^{-1}S \longrightarrow M_0 \setminus S$ is a homeomorphism;
\item if $p \in Z_H \subset Z_{(\lih)}$ and $q$ its image in
$S$ then $\pi^{-1}(q)$ is homeomorphic to $\pp W_p \git H^\cc$ where $W_p \leq T_pM$ is the Hermitian subspace defined in (\ref{w space});
\item $\widetilde M_0$ has a stratification indexed by $I_{\widetilde Z} = I_Z \setminus \{\lih\}$.
\end{enumerate}
\end{proposition}
Applying this proposition inductively we can find a symplectic
orbifold, arising as the regular reduction of a proper Hamiltonian $K$-space,
which we call the partial desingularisation of $M_0$. The reader is
referred to \cite{mes} for a precise uniqueness statement. As remarked
above, the only fact we will use is that it is unique up to homeomorphism. 
\section{The Kirwan map}
\label{kirwan map}
For a regular symplectic reduction $M_0$ of a proper Hamiltonian
$K$-space $M$ the Kirwan map is defined to be the
composition 
$$
\hk{*}{M} \longrightarrow \hk{*}{Z} \cong H^*(M_0)
$$
of restriction to the zero set $Z$ of the moment map and the natural
isomorphism of the equivariant cohomology of $Z$ with the cohomology
of $M_0$. One of its most important properties is that it is surjective, see
\cite[3.10]{k2}. In this section we generalise this and define a map 
$$
\hk{*}{M} \longrightarrow \ih{*}{M_0}
$$
to the intersection cohomology of the reduction even when the reduction
is not regular --- this idea was first discussed (for geometric
invariant theory quotients) in \cite{k3}. 

We introduce the notation $\ck{Z}$ for the derived push-forward of the constant sheaf with stalk $\qq$ via the map $EK \times_K Z \to M_0$. Our aim is to construct a morphism $\kappa : \ck{Z} \to \ic{M_0}$ in $\der{M_0}$. Inductively we may
assume there is a morphism 
$$
\widetilde \kappa : \ck{\widetilde Z} \longrightarrow \ic{\widetilde M_0}$$ 
where $\widetilde Z$ is the zero set of the moment map on $\widetilde M_0$. (In the base case of the induction when the reduction is regular there is a quasi-isomorphism.) Suppose that we have a `suitable' morphism $\rdf{\pi_*}\ic{\widetilde M_0} \to \ic{M_0}$.  We may then define $\kappa$ by the composition
$$
\ck{Z} \rightarrow \rdf{\pi_*}\ck{\widetilde Z}
\rightarrow \rdf{\pi_*}\ic{\widetilde M_0} \rightarrow \ic{M_0}
$$
of the equivariant pull-back $\pi^*$, the morphism $\rdf{\pi_*}\widetilde \kappa$ and this suitable morphism. We then define the Kirwan map to be the composite
$$
\hk{*}{M} \to \hk{*}{Z} \to \ih{*}{M_0}
$$
where the second map is induced from $\kappa$. In an abuse of
notation we will often denote this second map by $\kappa$ rather than the
correct but clumsy $\hyp{*}{M_0}{\kappa}$.

How do we conjure up this `suitable' morphism? 
Let us consider the inclusion $\pi^{-1}S \hookrightarrow
\widetilde M_0$ where $S$ was the maximal depth stratum of $M_0$ which we blew up in the first stage of the partial desingularisation. Here $\pi^{-1}S$ plays the role of exceptional
divisor in $\widetilde M_0$. In particular this inclusion is normally
nonsingular of codimension $2$. Thus, using the results of \cite[\S 5.4]{gm2}, we obtain a Gysin morphism
\begin{equation}
\label{gysin}
\rdf{\pi_*}\ic{\widetilde M_0} \stackrel{\gamma}{\longrightarrow} \rdf{\pi_*}\ic{\widetilde M_0}[2].
\end{equation}
\label{first chern class}
Intuitively we think of $\hyp{*}{M_0}{\gamma}$ as
analogous to multiplication by the first Chern class of the normal bundle. 

Recall from proposition \ref{pdprop} that if $q \in S$ is the image of $p
\in Z_H$ then $\pi^{-1}(q)$ is
homeomorphic to the geometric invariant theory quotient 
$\pp W_p \git H^\cc$.
Further the embedding of $\pi^{-1}(q)$ into the restriction of the normal bundle of $\pi^{-1}S \to \widetilde M_0$ to that fibre is homeomorphic to the embedding of the zero section into an anti-ample line bundle $\mathcal{L}$ on $\pp W_p \git H^\cc$.  
There is an isomorphism
$$
H_q^*(\rdf{\pi_*}\ic{\widetilde M_0}) \cong \ih{*}{\pp W_p \git H^\cc}. 
$$
Let $r = \frac{1}{2} \codim S = \frac{1}{2} \dim W_p \git H^\cc$. By \cite[\S 6]{bbd} there is a hard Lefschetz theorem for the above 
intersection cohomology groups arising from the ample line bundle
$\mathcal{L}^{-1}$. In particular there are isomorphisms  
$$
\ih{r -1-i}{\pp W_p \git H^\cc} \longrightarrow
\ih{r-1+i}{\pp W_p \git H^\cc} 
$$
for $i=0,\ldots,r-1$ given by $H^{r -1-i}_q(\gamma^i)$.
Moreover by considering the Gysin sequence for $\mathcal{L}$ we can identify the primitive part of $\ih{*}{\pp W_p \git H^\cc}$ as being isomorphic to $\ih{*}{W_p \git H^\cc}$. 

The decomposition theorem of \cite[\S 6]{bbd} can be thought of as a relative hard Lefschetz theorem. In this spirit we now show how the morphism $\gamma$ and the local Lefschetz decomposition above are sufficient to induce a global decomposition
\begin{equation}
\label{decomp}
\rdf{\pi_*}\ic{\widetilde M_0} \cong \ic{M_0} \oplus \dgs{B}
\end{equation}
where $\dgs{B}$ is supported on $S$. We will then take the projection to be the `suitable' morphism. For ease of reading we put $\dgs{A} = \rdf{\pi_*}\ic{\widetilde M_0}$ and set
$$
\dgss{A}{i} = \big(\tb{i}\rdf{\jmath_*}\jmath^*\dgs{A}\big)[2(i-r)] \qquad i=1,\ldots,r-1
$$
where $\jmath$ is the inclusion of $S$ in $M_0$.
Consider the morphism $\dgss{A}{1} \oplus \ldots \oplus \dgss{A}{r-1} \stackrel{\alpha}{\longrightarrow} \dgs{A}$ whose $i^{th}$ factor is 
$$
\dgss{A}{i} \to \rdf{\jmath_*}\jmath^*\dgs{A}[2(i-r)] \stackrel{\ \gamma^{r-i}}{\longrightarrow} \dgs{A}.
$$
For each $i <r-2$ we also have a morphism $(1,-\gamma):\dgss{A}{i}[2] \to \dgss{A}{i+1} \oplus \dgss{A}{i+2}$ and hence can define 
$$
\dgss{A}{1}[2] \oplus \ldots \oplus \dgss{A}{r-3}[2] \stackrel{\beta}{\longrightarrow} \dgss{A}{1} \oplus \ldots \oplus \dgss{A}{r-1}
$$
as the sum of these. It is clear that $\alpha\beta=0$ so we can choose a factorisation
\[
\xymatrix{
\dgss{A}{1} \oplus \ldots \oplus \dgss{A}{r-1} \ar[d] \ar[r]^{\qquad \quad \alpha} & \dgs{A}. \\
\con{\beta} \ar@{-->}[ur]_\varphi
}
\]
Since $\con{\beta}$ is supported on $S$, we get a morphism $\con{\varphi}\to \imath_*\imath^*\dgs{A}=\imath_*\imath^*\ic{M_0}$
where $\imath$ is the inclusion of the complement of $S$.
Computing cohomology we find that $\con{\varphi}\cong \tau_{<r}\con{\varphi}$ and that the induced morphism
$\con{\varphi}\cong \tau_{<r}\con{\varphi}\to \tau_{<r} \imath_*\imath^*\ic{M_0}\cong \ic{M_0}$ 
is a quasi-isomorphism. Since $\pi$ is  proper, Verdier duality applied to 
$$
\rdf{\pi_*}\ic{\widetilde M_0} \to \con{\varphi} \cong \ic{M_0}
$$ 
yields a morphism $\ic{M_0}\to\rdf{\pi_*}\ic{\widetilde M_0}$. Both are quasi-isomorphisms except on $S$ and it then follows from lemma \ref{unique morphisms} that they induce the desired splitting. Note that this decomposition is not canonical since we had to choose the lift $\varphi$.

Pictorially we think of this decomposition as follows (with r=4): 
\begin{center}
\begin{picture}(300,120)
\put(0,75){\line(1,0){70}}
\put(0,85){\line(1,0){10}}
\put(10,95){\line(1,0){10}}
\put(20,105){\line(1,0){10}}
\put(30,115){\line(1,0){10}}
\put(40,105){\line(1,0){10}}
\put(50,95){\line(1,0){10}}
\put(60,85){\line(1,0){10}}
\put(0,75){\line(0,1){10}}
\put(10,85){\line(0,1){10}}
\put(20,95){\line(0,1){10}}
\put(30,105){\line(0,1){10}}
\put(40,105){\line(0,1){10}}
\put(50,95){\line(0,1){10}}
\put(60,85){\line(0,1){10}}
\put(70,75){\line(0,1){10}}
\put(80,95){\makebox(0,0)[b]{$=$}}
\put(90,75){\line(1,0){20}}
\put(90,85){\line(1,0){10}}
\put(100,95){\line(1,0){10}}
\put(110,85){\line(1,0){10}}
\put(120,115){\line(1,0){10}}
\put(110,105){\line(1,0){10}}
\put(120,95){\line(1,0){10}}
\put(90,75){\line(0,1){10}}
\put(100,85){\line(0,1){10}}
\put(110,95){\line(0,1){10}}
\put(120,85){\line(0,1){10}}
\put(110,75){\line(0,1){10}}
\put(120,105){\line(0,1){10}}
\put(130,95){\line(0,1){20}}
\put(70,75){\line(0,1){10}}
\put(140,85){\makebox(0,0)[b]{$+$}}
\put(150,75){\line(1,0){50}}
\put(150,85){\line(1,0){10}}
\put(160,95){\line(1,0){10}}
\put(170,105){\line(1,0){10}}
\put(180,95){\line(1,0){10}}
\put(190,85){\line(1,0){10}}
\put(150,75){\line(0,1){10}}
\put(160,85){\line(0,1){10}}
\put(170,95){\line(0,1){10}}
\put(180,95){\line(0,1){10}}
\put(190,85){\line(0,1){10}}
\put(200,75){\line(0,1){10}}
\put(80,45){\makebox(0,0)[b]{$=$}}
\put(90,25){\line(1,0){20}}
\put(90,35){\line(1,0){10}}
\put(100,45){\line(1,0){10}}
\put(110,55){\line(1,0){10}}
\put(120,65){\line(1,0){10}}
\put(110,35){\line(1,0){10}}
\put(120,45){\line(1,0){10}}
\put(90,25){\line(0,1){10}}
\put(100,35){\line(0,1){10}}
\put(110,45){\line(0,1){10}}
\put(120,55){\line(0,1){10}}
\put(110,25){\line(0,1){10}}
\put(120,35){\line(0,1){10}}
\put(130,45){\line(0,1){20}}
\put(140,34){\makebox(0,0)[b]{$+$}}
\put(150,25){\line(1,0){30}}
\put(150,35){\line(1,0){10}}
\put(160,45){\line(1,0){10}}
\put(170,55){\line(1,0){10}}
\put(150,25){\line(0,1){10}}
\put(160,35){\line(0,1){10}}
\put(170,45){\line(0,1){10}}
\put(180,25){\line(0,1){30}}
\put(190,30){\makebox(0,0)[b]{$+$}}
\put(200,25){\line(1,0){20}}
\put(200,35){\line(1,0){10}}
\put(210,45){\line(1,0){10}}
\put(200,25){\line(0,1){10}}
\put(210,35){\line(0,1){10}}
\put(220,25){\line(0,1){20}}
\put(230,28){\makebox(0,0)[b]{$+$}}
\put(240,25){\line(1,0){10}}
\put(240,35){\line(1,0){10}}
\put(240,25){\line(0,1){10}}
\put(250,25){\line(0,1){10}}
\put(260,27){\makebox(0,0)[b]{$-$}}
\put(270,25){\line(1,0){10}}
\put(270,35){\line(1,0){10}}
\put(270,25){\line(0,1){10}}
\put(280,25){\line(0,1){10}}
%
%
\linethickness{0.001mm}
\multiput(272,25)(2,0){5}{\line(0,1){10}}
\multiput(270,27)(0,2){5}{\line(1,0){10}}
\multiput(202,25)(2,0){5}{\line(0,1){10}}
\multiput(200,27)(0,2){5}{\line(1,0){10}}
\multiput(170,25)(2,0){5}{\line(0,1){10}}
\multiput(170,27)(0,2){5}{\line(1,0){10}}
\put(140,0){\makebox(0,0)[b]{``\ $\rdf{\pi_*}\ic{\widetilde M_0}|_S \cong \ic{M_0}|_S+\gamma \dgss{A}{3}+\gamma^2 \dgss{A}{2}+\gamma^3 \dgss{A}{1} - \gamma^2 \dgss{A}{1}$\ ''}}
\end{picture}
\end{center}

\begin{conjecture}
The Kirwan map $\hk{*}{M} \to \ih{*}{M_0}$ is surjective.
\end{conjecture}
There is various evidence suggesting that this conjecture should be
true. First note that it follows from the equivariant Morse theory of
the moment map that $\hk{*}{M} \rightarrow \hk{*}{Z}$ is surjective (see
\cite[3.10]{k1}) so that what is in question is the
surjectivity of $\hk{*}{Z} \rightarrow \ih{*}{M_0}$. In various special
cases, for instance when the action of $K$ is almost-balanced or $K$
is a circle, we show below that this is surjective. 
In the GIT setting, the second author proved the surjectivity \cite{woolf}
by using the decomposition theorem.
Recent work of Tolman suggests that this is so for arbitrary actions by
higher dimensional tori.
\section{The cosupport condition}
\label{cosupport condition}
When $K$ acts quasi-freely on $Z$ there is a quasi-isomorphism $\ic{M_0} \cong \ck{Z}$. In other words $\ck{Z}$ satisfies both the support and cosupport conditions. Since the equivariant cohomology of a point with respect to a non-trivial connected group is infinite dimensional we can easily check that
$$
\ck{Z} \textrm{ satisfies support } \iff K \textrm{ acts quasi-freely on } Z.
$$
It is very natural to ask: when does $\ck{Z}$ satisfy cosupport? Again there is a pleasant geometric interpretation of this condition. 
\begin{definition}
For each connected component of a stratum of $M_0$ fix a point with generic stabiliser in the preimage of the component in $Z$. We say that the \emph{action satisfies $C(j)$} if, for each such point $p$,
$$
\dim(W_p \git H^\cc) + j< 2 \min\{\codim_{\pp W_p} S \ | \ S \in \mathcal{U}_p\},
$$
where $H = \Stab p$, the subspace $W_p$ of $T_pM$ is as in (\ref{w space}) and $\mathcal{U}_p$ is the set of unstable strata in the Morse stratification of $\pp W_p$ by the gradient flow of the norm square of the moment map associated to the $H$ action. Equivalently, using the connection with geometric invariant theory outlined in \S\ref{git}, we can express this as
$$
\dim(W_p \git H^\cc)  + j <  2 \codim_{W_p}\{x \in W_p \ |\  0 \in \overline{H^{\cc}x}\}.
$$
Note that $\{x \in W_p \ |\  0 \in \overline{H^{\cc}x}\}$ is $\varphi_p^{-1}(0)$ where $\varphi_p : W_p \to W_p \git H^\cc$ is the (algebraic) quotient map. This definition is independent of the chosen points in the preimages of the connected components of the strata.  
\end{definition}
\begin{lemma}
$\ck{Z} \in \pdera{m^+}{j}{M_0} \iff \textrm{ the action satisfies } C(j).
$
\end{lemma}
\begin{proof}
Suppose $p$ is in the preimage in $Z$ of a point $q$ in a connected component $\jmath_S : S \hookrightarrow M_0$ of a stratum, and that $H=\Stab p$. There is a long exact sequence
$$
\ldots \to H_q^i(\jmath_S^! \ck{Z}) \to H^i_H(W_p) \to H^i_H(W_p \setminus \varphi_p^{-1}(0)) \to \ldots
$$
We see that $\codim \varphi_p^{-1}(0) = \min\{i | H_q^i(\jmath_S^! \ck{Z})\neq 0\}$. By the definition of $\pdera{m^+}{j}{M_0}$ we see that $\ck{Z}$ lies in this subcategory if, and only if, 
$$
\codim \varphi_p^{-1}(0) > \frac{1}{2}\codim S + j = \frac{1}{2}\dim(W_p\git H^\cc) + j
$$
which is precisely $C(j)$.
\end{proof}
In particular, since $\ck{Z}$ satisfies cosupport if, and only if, it lies in the subcategory $\pdera{m^+}{0}{M_0}$, we have 
$$
\ck{Z} \textrm{ satisfies cosupport } \iff \textrm{ the action satisfies } C(0).
$$
\begin{lemma}$C(0) \implies \ptb{m^-}{0}\ck{Z} \cong \ic{M_0} \implies C(-1).$
\end{lemma}
\begin{proof}
This follows from the triangle associated to $\ptb{m^-}{0}\ck{Z} \to \ck{Z}$ and the above lemma.
\end{proof} 
\begin{remark} If the condition $C(0)$ is satisfied, we say the $K$ action
on $M$ is \emph{almost-balanced}. The terminology came from the fact that
$C(0)$ can be viewed as a condition on the distribution of
the weights of the maximal torus action on the normal space to each stratum \cite{kiem}.
For example, if the weights of the maximal torus action of $H$ on $W_p$ is symmetric with respect to the origin, for each $\lih\in I_Z$ and $p\in Z_{\lih}$,
then the $K$ action is almost balanced. This is the case for
the quotient construction yielding the moduli space of bundles on a Riemann surface (\cite{kiem} Proposition 7.3). We compute a simpler example in \S \ref{simple example}.  
\end{remark}
When the action is almost-balanced the above lemma shows that there is a natural morphism $\lambda : \ic{M_0} \cong \ptb{m^-}{0}\ck{Z} \to \ck{Z}$
which we interpret as a generalisation of the pull-back which exists when the  action of $K$ on $Z$ is quasi-free. 

\begin{proposition}
\label{splitting}
Suppose the action of $K$ on $M$ is almost balanced. Then $\kappa\lambda = \id$ and so the intersection cohomology of $M_0$ is naturally identified with a subspace of the equivariant cohomology of $Z$.
\end{proposition}
\begin{proof}
It is easy to check that the composition 
$$
\ic{M_0} \stackrel{\lambda}{\longrightarrow} \ck{Z} \stackrel{\kappa}{\longrightarrow} \ic{M_0}
$$
is non-zero, but then it must then be a quasi-isomorphism by \ref{unique morphisms}.
\end{proof}
\subsection{Finding the right subspace}
\label{finding}
In order to make practical use of the observation that $\ih{*}{M_0}$ embeds into $\hk{*}{Z}$ when the action is almost-balanced we need to be able to give a concrete description of the image. To this end we introduce an auxiliary complex $\ve{Z}$.

Suppose $\lih$ is in the indexing set $I_Z$ for the stratification of $M_0$. Let $Y_{(\lih)}$ be the closure in $Z$ of the corresponding subset $Z_{(\lih)}$, and similarly let $Y_\lih$ be the closure of $Z_\lih$. Let $N^{H_0}$ be the normaliser in $K$ of $H_0=\mathrm{exp}\, \lih$, the identity component of the generic stabiliser of points in $Z_\lih$. Let $N^{H_0}_0$ be the identity component of $N^{H_0}$.

Notice that 
$Y_{(\lih)}$ is homeomorphic to $K Y_\lih$ and consider the resolution 
$$K\times_{N^{H_0}}Y_{\lih}\to KY_{\lih}=Y_{(\lih)}$$
which is an isomorphism over $KZ_{\lih}$.
Since the normal subgroup $H_0$ of $N^{H_0}_0$ acts trivially upon $Y_\lih$, we get morphisms 
\begin{equation}
\label{qis}
\ck{Y_{(\lih)}} \to \eqc{N^{H_0}}{Y_\lih} \cong \rdf{{\psi_\lih}_*} \big(\eqc{N^{H_0}_0/H_0}{Y_\lih} \otimes H\udot_{\!\!\!H_0}\big)
\end{equation}
where $H\udot_{\!\!\!H_0}$ is the constant sheaf on $Y_\lih/N^{H_0}_0$ with stalk $H^*_{H_0}$ and $\psi_\lih$ is the quotient map $Y_\lih / N^{H_0}_0 \to Y_\lih / N^{H_0}$. Define the complex $\mathcal{L}_\lih(Z) \in \der{M_0}$ to be 
the extension by zero to $M_0$ of 
$$
\rdf{{\psi_\lih}_*} \big(\eqc{N^{H_0}_0/H_0}{Y_\lih} \otimes \tgeq{n_\lih}H\udot_{\!\!\!H_0}\big)
$$
where $n_\lih$ is half the real codimension of the stratum $Z_{(\lih)}/K$ of $M_0$.
\begin{definition}
\label{defn of v}
Define $\ve{Z}$ (up to quasi-isomorphism) by the distinguished triangle
$$
\ve{Z} \to \ck{Z} \to \bigoplus_{\lih \in I_Z} \mathcal{L}_\lih(Z)
$$
where the morphism $\ck{Z} \to \mathcal{L}_\lih(Z)$ is given by restricting to $Y_{(\lih)}$, applying the morphisms in (\ref{qis}) and truncating in the second factor. 

Let $V^*(Z)$ be the image of the hypercohomology of $\ve{Z}$ in $H^*_K(Z)$. Restriction to $Y_{(\lih)}$ composed with the morphisms (\ref{qis}) induces a map
$$
H^*_K(Z) \to [H^*_{N^{H_0}_0/H_0}(Y_\lih) \otimes H^*_{H_0}]^{\pi_0 N^{H_0}}
$$
where the square brackets denote the invariant part under the action of the finite group $\pi_0 N^{H_0} \cong N^{H_0} / N^{H_0}_0$. We can check that $V^*(Z)$ consists of those classes in $H^*_K(Z)$ whose image under this map lies in the direct summand
$$
[H^*_{N^{H_0}_0/H_0}(Y_\lih) \otimes H^{<n_\lih}_{H_0}]^{\pi_0 N^{H_0}}
$$
for each $\lih \in I_Z$.
\end{definition}
\begin{assumption}
\label{wb assumption}
For the remainder of this section we will make the technical assumption that the action is not only almost-balanced 
but in fact \emph{weakly-balanced} i.e. not only does $\ck{Z}$ satisfy cosupport (which means $\ck{Z}\in {\bf D} ^{m+}_{\ge 0}(M_0)$ )
but so does each of the $\eqc{N^{H_0}_0/H_0}{Y_\lih}$ (by which we mean $\eqc{N^{H_0}_0/H_0}{Y_\lih}\in {\bf D} ^{m+}_{\ge 0}(Y_{\lih}/N^{H_0})$,
or equivalently 
$\mathcal{L}_\lih(Z)\in {\bf D} ^{m}_{\ge 1}(M_0)$). These latter conditions can be interpreted geometrically in a similar fashion to the interpretation of the almost-balanced condition --- see \cite{kiem}. The moduli space of bundles on a Riemann surface arises as a quotient by a weakly-balanced action and we compute another example in \S \ref{simple example}.
\end{assumption}
\begin{lemma}
\label{v splits}
Suppose the action of $K$ is weakly-balanced. Then there is a direct sum decomposition $\ve{Z} \cong \ic{M_0} \oplus \dgsx{F}{Z}$.
\end{lemma}
\begin{proof}
It follows immediately from the weakly-balanced assumption that the morphism $\ic{M_0} \cong \ptb{m^-}{0}\ck{Z} \to \ck{Z}$ can be factored through $\ve{Z}$ since $\tau ^{m-}_{\le 0}\mathcal{L}_\lih(Z)=0$. It is easy to check that the resulting composition
$$
\ic{M_0} \to \ve{Z} \to \ck{Z} \stackrel{\kappa}{\longrightarrow} \ic{M_0} 
$$
is not zero and so must then be a quasi-isomorphism by lemma \ref{unique morphisms}.
\end{proof}
Thus the intersection cohomology $\ih{*}{M_0}$ embeds as a subspace of $V^*(Z) \subset \hk{*}{Z}$. We now show that the image of the embedding is precisely this subspace.
\begin{theorem}
\label{ih is v}
Suppose the action of $K$ is weakly-balanced. Then the restriction of the Kirwan map to $V^*(Z)$ is an isomorphism onto $\ih{*}{M_0}$.
\end{theorem}
\subsection{Proof of theorem \ref{ih is v}}
Let $F^*(Z)$ be the image of the hypercohomology of $\dgsx{F}{Z}$ in $\hk{*}{Z}$. It follows from lemma \ref{v splits} that it is sufficient to show that $F^*(Z)$ is zero. Inductively we will assume that $F^*(\widetilde Z)$ is zero, where 
\begin{equation}
\label{meq1}
\ve{\widetilde Z} \cong \ic{\widetilde M_0} \oplus \dgsx{F}{\widetilde Z}
\end{equation}
and $F^*(\widetilde Z)$ is the image of the hypercohomology of $\dgsx{F}{\widetilde Z}$ in $\hk{*}{\widetilde Z}$.

Let $S$ be the deepest stratum of $M_0$ (which is the centre of the first blowup in the partial desingularisation procedure) and set $r= \frac{1}{2}\codim S$. We also introduce closed subsets $C \subset Z\label{c and e}$ (for centre) and $E \subset \widetilde Z$ (for exceptional divisor) which are respectively the preimages of $S$ under $Z \to M_0$ and $\widetilde Z \to \widetilde M_0 \to Z$. 
\begin{lemma}
\label{lem1}
We can choose a factorisation
\[
\xymatrix{
\dgsx{V}{Z} \ar@{-->}[r]^\alpha \ar[dr] & \tbon{r}{S}\ck{Z} \ar[d]\\
& \ck{Z}
}
\]
where $\tbon{r}{S}$ is the truncation over a closed subset functor (see \ref{truncation}).
\end{lemma}
\begin{proof}
This follows immediately from definition \ref{defn of v} of $\dgsx{V}{Z}$. 
\end{proof}
\begin{lemma}
\label{lem2}
We can choose a morphism $\beta$ so that
\[
\xymatrix{
\dgsx{V}{Z} \ar@{-->}[r]^\beta \ar[d] & \rdf{\pi_*}\ve{\widetilde Z} \ar[d] \\
\ck{Z} \ar[r] & \rdf{\pi_*} \ck{\widetilde Z}
}
\]
commutes. 
\end{lemma}
\begin{proof}
For each $\lih \in I_{\widetilde{Z}}$, let $\widetilde{Y}_{\lih}$ denote the $H_0=\mathrm{exp}\,\lih$ fixed point set in $\widetilde{Z}$ just like $Y_{\lih}$ in $Z$.
Then since the map $\widetilde{Z}\to Z$ is equivariant, we have the commutative diagram
$$\xymatrix{
K\times_{N^{H_0}}\widetilde{Y}_{\lih} \ar[r]\ar[d] & K\widetilde{Y}_{\lih}\ar[d]\ar@{^(->}[r] & {\widetilde{Z}}\ar[d]\\
K\times_{N^{H_0}}{Y}_{\lih} \ar[r] & K{Y}_{\lih}\ar@{^(->}[r] & {Z}.
}$$
This induces the following commutative diagram
$$\xymatrix{
{\ck{Z}}\ar[r]\ar[d] & {\rdf{\pi_*}\ck{\widetilde{Z}}}\ar[d]\\
{\bigoplus_{\lih\in I_Z}\mathcal{L}_{\lih}(Z)}\ar[r] & {\bigoplus_{\lih\in I_{\widetilde{Z}}} \rdf{\pi_*}\mathcal{L}_{\lih}(\widetilde{Z})}
}$$
by assigning zero for $\lih\in I_Z\setminus I_{\widetilde{Z}}$. The claim now follows from the definitions.
\end{proof}
\begin{lemma}
\label{lem3}
If $\eta \in F^*(Z)$ then $0 = \pi^*\eta|_E \in \hk{*}{E}$ where $E$ is defined above.
\end{lemma}
\begin{proof}
We can decompose $\rdf{\pi_*}\ve{\widetilde Z}$ into a direct sum $\rdf{\pi_*}\ic{\widetilde M_0} \oplus \rdf{\pi_*} \dgsx{F}{\widetilde Z}$ using (\ref{meq1}). It follows from lemma \ref{v splits} that the composition 
$$
\rdf{\pi_*}\ic{\widetilde M_0} \oplus \rdf{\pi_*} \dgsx{F}{\widetilde Z} \cong \rdf{\pi_*}\ve{\widetilde Z} \to \rdf{\pi_*}\ck{\widetilde Z} \to \rdf{\pi_*}\ic{\widetilde M_0}
$$
is a quasi-isomorphism on the first term. Now let us consider the composition
$$
\dgsx{F}{Z} \to \dgsx{V}{Z} \stackrel{\beta}{\longrightarrow} \rdf{\pi_*}\dgsx{V}{\widetilde Z} \cong  \rdf{\pi_*}\ic{\widetilde M_0} \oplus \rdf{\pi_*}\dgsx{F}{\widetilde Z}.
$$
By the above and lemma \ref{lem2} the projection onto the first term of the RHS is, up to a quasi-isomorphism, the same as 
$$
\dgsx{F}{Z} \to \dgsx{V}{Z} \to \ck{Z} \to \rdf{\pi_*}\ck{\widetilde Z} \to \rdf{\pi_*}\ic{\widetilde M_0}.
$$
In particular if we decompose $\rdf{\pi_*}\ic{\widetilde M_0}$ as in (\ref{decomp}) into $\ic{M_0} \oplus \dgs{B}$ then by definition $\dgsx{F}{Z} \to \ic{M_0}$ is zero and the choice of $\alpha$ in lemma \ref{lem1} induces a unique lift of $\dgsx{F}{Z} \to \dgs{B}$ to $\tbon{r}{S}\dgs{B}$. 

The inductive assumption that $F^*(\widetilde Z)=0$ now tells us that the two morphisms
\[
\xymatrix{
&  \tbon{r}{S} \ck{Z}\ar[dr] && \\
\dgsx{F}{Z} \ar[ur]^{\alpha'} \ar[dr]_{\beta'} && \tbon{r}{S}\rdf{\pi_*}\ck{\widetilde Z} \ar[r] & \rdf{\pi_*}\ck{\widetilde Z} \\
& \tbon{r}{S}\dgs{B} \ar[ur]
}
\]
induced by $\alpha$ and $\beta$ respectively must give the same map on hypercohomology.

The almost-balanced assumption guarantees that there is a quasi-isomorphism $\tbon{r}{S} \ck{Z} \cong \tbon{r}{S}\rdf{\imath_*}\imath^*\ck{Z}$, where $\imath$ is the inclusion of the open complement to $S$. So $\tbon{r}{S}\rdf{\pi_*}\ck{\widetilde Z}$ decomposes as $\tbon{r}{S} \ck{Z} \oplus \dgs{C}$ where $\dgs{C}$ is defined (up to quasi-isomorphism) by the distinguished triangle
$$
\dgs{C} \to \tbon{r}{S}\rdf{\pi_*}\ck{\widetilde Z} \to  \tbon{r}{S}\rdf{\imath_*}\imath^*\ck{Z}.
$$
Furthermore because $\dgs{B}$ is supported on $S$ it includes into $\tbon{r}{S}\rdf{\pi_*}\ck{\widetilde Z}$ as a subobject of $\dgs{C}$. 

The final observation we require is that the spectral sequence computing the hypercohomology of $\rdf{\pi_*}\ck{\widetilde Z}|_S$ degenerates at the $E_2$-term. This follows from Deligne's criterion and the surjectivity  of the restriction $H^*_H(\pp W_p)\to H^*_H(\pp W_p^{ss})$. Hence the spectral sequence for $\tbon{r}{S} \rdf{\pi_*}\ck{\widetilde Z}|_S$ must also degenerate and there is an injection 
\begin{equation}
\label{meq2}
\hyp{*}{S}{\tbon{r}{S} \rdf{\pi_*}\ck{\widetilde Z}|_S} \hookrightarrow 
\hyp{*}{S}{ \rdf{\pi_*}\ck{\widetilde Z}|_S} \cong \hk{*}{E}
\end{equation}
where $E$ is the `exceptional divisor' in $\widetilde Z$. But now we are done for we have shown that the image of the class $\eta \in F^*(Z)$ pulled back to $\hk{*}{E}$ lies in two subspaces, namely $\hyp{*}{S}{\tbon{r}{S} \ck{Z}|_S}$ and $\hyp{*}{S}{\dgs{C}}$, whose intersection is zero.
\end{proof}
We can now complete the proof of theorem \ref{ih is v}. Let $\eta \in F^*(Z)$. By lemmas \ref{lem2} and \ref{lem3} we see that $\pi^*\eta \in \hk{*}{\widetilde Z}$ lies in the subspace $V^*(\widetilde Z) \cong \ih{*}{\widetilde M_0}$ and restricts to zero on $E$. Furthermore by the definition of $F^*(Z)$ we know that it projects to zero in the summand $\ih{*}{M_0}$ of $\ih{*}{\widetilde M_0}$. It follows from (\ref{decomp}) and the fact that $\dgs{B}$ is supported on $S$ that we must have $\pi^*\eta = 0$. 

Studying the proof of lemma \ref{lem3} more carefully we see that we can actually conclude that not only is $\pi^*\eta|_E=0$ but that $\eta|_C = 0$ where $C\subset Z$ is the `centre' of the blowup defined above on page \pageref{c and e} since it lies in $\hyp{*}{S}{\tbon{r}{S} \ck{Z}|_S}$ which injects into
$H^*_K(E)$. 
Finally we note that 
\begin{equation}\label{lastint}
\ker \big(\hk{*}{Z} \to \hk{*}{\widetilde Z} \big) \cap \ker \big(\hk{*}{Z} \to \hk{*}{C}\big) = 0
\end{equation}
so that $\eta = 0$ as desired. (\ref{lastint}) follows from the observation
that if $\tilde{M}\to U$ is the blowup of a neighborhood of $Z$ along
a submanifold containing $C$ such that $\tilde{M}\git K=\tilde{M}_0$,
then the nonminimal strata in $\tilde{M}$ with respect to the norm square
of the moment map retract onto those in the exceptional divisor $\tilde{E}$
and hence $\ker \big( \hk{*}{\tilde{M}}\to\hk{*}{\widetilde Z} \big)
\cong \ker \big( \hk{*}{\tilde{E}}\to\hk{*}{E} \big)$.
\subsection{The intersection pairing}
\label{intersection pairing}
Since $M_0$ is compact (by the properness assumption on the moment map) there is an intersection pairing
$$
\ih{*}{M_0} \otimes \ih{*}{M_0} \longrightarrow \qq.
$$
This arises from a morphism $\ic{M_0}^{\otimes 2}\to \ict{M_0}$ where the latter is the intersection cohomology complex with top perversity $t(S)=\codim(S)$. When the action is almost-balanced there is a second natural morphism 
$$
\ic{M_0}^{\otimes 2}\to \ck{Z}^{\otimes 2} \to \ck{Z} \to \ic{M_0} \to \ict{M_0}
$$
induced from the product on equivariant cohomology. 
\begin{lemma}
\label{morphisms agree}
These two morphisms $\ic{M_0}^{\otimes 2}\to\ict{M_0}$ are the same.
\end{lemma}
\begin{proof}
It is easy to check that $\ic{M_0}^{\otimes 2} \in \pderb{t^-}{0}{M_0}$. Then by lemma \ref{unique morphisms} we have
$$
\Hom(\ic{M_0}^{\otimes 2},\ict{M_0}) \cong \qq.
$$
Since on the non-singular stratum the two morphisms are both just the product on the constant sheaf with stalk $\qq$ we are done.
\end{proof}
\begin{lemma} 
\label{compactly supported}
Let $\Sigma$ be the union of the singular strata of $M_0$. Then 
$$
\hyp{\dim M_0}{\Sigma}{\ic{M_0}^{\otimes 2}}=0.
$$
\end{lemma}
\begin{proof}
Note that the maximal degree in which the hypercohomology of a complex $\dgs{C}$ does not vanish is bounded by 
$$
\max_i\{i+\dim \supp \sh{i}{\dgs{C}}\}.
$$
Now since $\ic{M_0}^{\otimes 2} \in \pderb{t^-}{0}{M_0}$ we see that this quantity is strictly less than $\max\{\dim S_\alpha + t(S_\alpha)\} = \dim M_0$.
\end{proof}
\begin{remark}
\label{pairings remark}
Let $d=\dim M_0$. It follows from the above lemma that the product on equivariant cohomology gives us maps
$V^i(Z) \otimes V^{d -i}(Z) \to V^d(Z)$
and, by a similar argument, that $V^d$ is the image of the pull-back $\h{d}{M_0} \to \hk{d}{Z}$. 
\end{remark}
\begin{definition}
We say a coadjoint orbit $\mathcal{O}$ is \emph{close to $0$} if $\mu^{-1}\mathcal{O}$ is contained in the set $M^{ss}$ of points which flow to $Z$ under the gradient flow of the norm square of the moment map. (This depends upon the choices of compatible metric on $M$ and invariant metric on $\lik$.)
\end{definition}
\begin{proposition}
Suppose the coadjoint orbit $\mathcal{O}$ is close to $0$ and contains a regular value of $\mu$ (and hence consists entirely of regular values). Let $\eta,\zeta \in \hk{*}{M}$ be two classes whose restrictions to $\hk{*}{Z}$ lie in $V^*(Z)$. Then the intersection pairing of $\kappa(\eta|_Z)$ with $\kappa(\zeta|_Z)$ in $\ih{*}{M_0}$ is given by 
$$
\eta\zeta|_{\mu^{-1}(\mathcal{O})} [M_\mathcal{O}].
$$
\end{proposition}
\begin{proof}
Since $\mu^{-1}(\mathcal{O}) \subset M^{ss}$ the flow induces an equivariant map $\mu^{-1}(\mathcal{O}) \to \mu^{-1}(0) =Z$ which is a homeomorphism between dense open sets and so we have
\[
\xymatrix{
\hk{*}{Z} \ar[d] & \h{*}{M_0} \ar[l] \ar[d]\\
\hk{*}{\mu^{-1}(\mathcal{O})} & \h{*}{\mu^{-1}(\mathcal{O})/K}. \ar[l]
}
\] 
Let $d = \dim M_0$. By lemma \ref{morphisms agree} and remark \ref{pairings remark} we also have a commutative diagram
\[
\xymatrix{
V^i(Z) \otimes V^{d -i}(Z) \ar[d] & \ih{i}{M_0} \otimes \ih{d-i}{M_0} \ar[d] \ar[l] \\
V^d(Z) & \h{d}{M_0} \ar[l]
}
\]
where the horizontal maps are isomorphisms. The result follows by composing these two diagrams and noting that $\h{d}{M_0} \to \h{d}{\mu^{-1}(\mathcal{O})/K}$ preserves evaluations against the fundamental class because $\mu^{-1}(\mathcal{O})/K \to M_0$  is `birational'.
\end{proof} 
\subsection{Example}
\label{simple example}
({cf. \cite[\S 4]{k3}})
Let $M = (\pp^1)^n$ acted on diagonally by $SU_2$ (in the obvious way). Identifying $\pp^1$ with $S^2 \subset \rr^3$ and $\mathfrak{su}_2$ with $\rr^3$ we have a moment map
\begin{eqnarray*}
\mu : M & \longrightarrow & \rr^3 \\
(x_1, \ldots , x_n) & \longmapsto & \sum x_i
\end{eqnarray*}
for this action. We can check that this action is almost-balanced. The zero set $Z = \mu^{-1}(0)$ consists of $n$-tuples of points balanced about the origin. Let $T \subset SU_2$ be the maximal torus
$$
\left( 
\begin{array}{cc} 
e^{it} & 0\\
0 & e^{-it}
\end{array}
\right).
$$
The Weyl group of $SU_2$ is $\zz_2$ and $H^*_{SU_2}(M) = [H^*_T(M)]^{\zz_2}$
is the $\zz_2$-invariant part of $H^*_T(M)$. The $T$-equivariant cohomology is generated by classes $\alpha_1,\ldots,\alpha_n,\beta$ in degree $2$ subject to the relations $\alpha_i^2 = \beta^2$. The Weyl group acts trivially on the $\alpha_i$ and takes $\beta$ to $-\beta$. Hence 
$$
H^*_{SU_2}(M) \cong \cc[\alpha_1,\ldots, \alpha_n,\beta^2] / \langle \alpha_i^2 = \beta^2 : i=1,\ldots, n \rangle.
$$
The only non-discrete subgroups which have fixed point are conjugates of $T$. The fixed points of $T$ correspond to sequences $(a_1,\ldots,a_n)$ where $a_i= \pm 1$ respectively according to whether we take $0$ or $\infty$ in the $i^{th}$ $\pp^1$. If $F$ is the fixed point corresponding to $(a_1,\ldots,a_n)$ then $\mu_T(F) = \sum a_i$, the $T$-equivariant Euler class $e_F$ of the normal bundle to $F$ is $(\prod a_i)\beta^n$ and $\alpha_i |_F = a_i \beta \in H^*_T(F)$. 

Since $\mu(F)=0$ implies that $\mu_T(F)=0$ it is clear that the reduction at zero is seriously singular (i.e. has worse than orbifold singularities) if, and only if, $n$ is even. Put $n=2m$. We note that the singularities of the reduction $M_0$ of $M$ are isolated. So the restriction of $\eta \in H^*_{SU_2}(M)$ to $Z$ lies in $V^*$ if, and only if, either $\deg \eta < 2m-3$ or $\eta|_{Z_{(\lit)}}=0$. Furthermore it follows from the Poincar\'e polynomial calculations of \cite[\S 5]{k1} that the restriction
$$
H^r_{SU_2}(M) \longrightarrow H^r_{SU_2}(Z)
$$
is an isomorphism for $r<2m$. We deduce that, for $r<2m-3$, there is an isomorphism $V^r \cong H^r_{SU_2}(M)$. 

The reduction $M_0$ can, via geometric invariant theory, be given a description as a projective algebraic variety. In particular its intersection cohomology will satisfy the hard Lefschetz theorem i.e. there are isomorphisms
$$
\ih{2m-3-i}{M_0} \stackrel{\sim}{\longrightarrow} \ih{2m-3+i}{M_0}
$$
given by powers of a map $\ih{r}{M_0} \to \ih{r+2}{M_0}$. A moment's thought shows that on $V^*$ this map must be given  by multiplication by the K\"ahler class $\omega = \sum \alpha_j$. (Note that any multiple of $\omega$ restricts to an element of $V^*$ since 
$$
\omega |_{Z_{(\lit)}} = \sum_{\sum a_j = 0} \sum a_j \beta  = 0.)
$$
Finally since the equivariant cohomology of $M$ vanishes in odd degree we deduce that so does that of $M_0$ and, in particular, $\ih{2m-3}{M_0}=0$. Hence
$$
\ih{r}{M_0} \cong \left\{
\begin{array}{ll}
H^r_{SU_2}(M) & r = 2m-3-i , \ i\geq 0 \\
\omega^i H^{r-2i}_{SU_2}(M) & r=2m-3+i, \ i\geq 0
\end{array}
\right.
$$
with the intersection pairing given by taking the product in $H^*_{SU_2}(M)$ and evaluating against the fundamental class of the reduction of $M$ at any regular coadjoint orbit close to $0$. The full Lefschetz decomposition can be computed by considering the action of powers of $\omega$. 
\section{Circle actions}
\label{circle actions}
As we would expect the circle case is far simpler than the general one. Let $M$ be a proper Hamiltonian $S^1$-space, $Z$ be the zero set of an equivariant moment map and $M_0=Z/S^1$ be the reduction. Let $\mathcal{F}_0$ be the set of fixed point components of $S^1$ which lie in $Z$. Each such component corresponds to an isomorphic singularity in $M_0$. For each $F \in \mathcal{F}_0$ define $d(F)$ (respectively $e(F)$) to be the minimum (respectively maximum) of the numbers of positive and negative weights of $S^1$ on  normal fibre $W_F$ to $F$ in $M$. 

The normal fibre to the singularity corresponding to $F \in \mathcal{F}_0$ in $M_0$ is given by the quotient $W_F \git \cc^*$. We can easily check that removing the vertex we have
$$
(W_F\setminus \{0\})\git \cc^* \cong (W_F^+\setminus\{0\}) \times_{\cc^*} (W_F^-\setminus\{0\})
$$ 
where $W_F^\pm$ are the positive and negative weight spaces. This is the total space of some $\cc^{e(F)} \setminus \{0\}$ bundle over a weighted projective space of dimension $d(F)-1$. 

Define a perversity $n$ by 
$$
n(S) = \left\{
\begin{array}{ll}
0 & \codim S = 0 \\
2d(F)-1 & S \cong F.
\end{array}
\right.
$$
\begin{lemma}
There is a quasi-isomorphism $\ic{M_0} \cong \icp{n}{M_0}$ and $\ck{Z}$ lies in the full subcategory $\pdera{n^+}{0}{M_0}$.
\end{lemma}
\begin{proof}
Both statements follow almost immediately from the above calculation of the normal fibre.
\end{proof}
We now have the same situation as in the last section but with the perversity $n \leq m$ replacing the middle perversity $m$. Arguing analogously we deduce
\begin{theorem}
The intersection cohomology $\ih{*}{M_0}$ is isomorphic (via the Kirwan map) to the subspace
$$
\{ \eta \in H^*_{S^1}(Z) \ | \ \eta|_F \in H^*(F) \otimes H^{<2d(F)-1}_{S^1}\}.
$$
The intersection pairing is given by taking the product of classes and evaluating against the fundamental class of either of the shift resolutions $M_{\pm\epsilon}$ where $0 \neq \epsilon \in \rr$. 
\end{theorem}
\begin{corollary}
The intersection Betti numbers of $M_0$ are given by the Poincar\'e polynomial
$$
P^{S^1}_t(Z) - \frac{1}{1-t^2}\sum_{F\in\mathcal{F}_0}t^{2d(F)}P_t(F).
$$
\end{corollary}
\begin{example}
We consider a linear circle action on the projective space $\pp^n$.
Let $p, q$ and $r$ denote the number of positive, negative, and zero weights respectively,
so that $n=p+q+r-1$. Let us assume that $p\leq q$ (the other case being entirely similar).
Using equivariant Morse theory, we get
$$
\begin{array}{rl}
P^{S^1}_t(Z)&=(P_t(\pp^n)-t^{2q+2r}P_t(\pp^{p-1})
-t^{2p+2r}P_t(\pp^{q-1}))/(1-t^2)\\
&=(1+t^2+\cdots +t^{2p+2r-2}-t^{2q+2r}-\cdots -t^{2n})/(1-t^2).
\end{array}
$$
Hence, by the above corollary, the intersection Poincar\'e polynomial is
$$
P^{S^1}_t(Z)-\frac{t^{2p}P_t(\pp^{r-1})}{1-t^2} =\frac{(1-t^{2p})(1-t^{2q+2r})}{(1-t^2)^2}
$$
which is a palindromic polynomial of degree $2n-2$.

Intersection pairings can be also computed by using the splitting above ---
see \cite{kiem}.
\end{example}

\end{document}